\documentclass[twoside,12pt]{article}
\usepackage{amssymb,amsmath,euscript}
\usepackage{graphicx,color,caption}
\usepackage{subfigure}
\usepackage{amsmath,amsfonts,amssymb,amsthm,graphicx,booktabs,enumerate,mathrsfs,color,subfigure,multicol,mathdots}
\usepackage{algorithm}
\usepackage{epstopdf}
\usepackage{algorithmic}
\usepackage{amssymb,amsmath,euscript}
\usepackage{graphicx,color,caption}
\usepackage[all]{xy}
\usepackage{cite}
\newtheorem{proposition}{Proposition}[section]
\newtheorem{theorem}[proposition]{Theorem}
\newtheorem{lemma}[proposition]{Lemma}

\newtheorem{definition}[proposition]{Definition}

\newtheorem{example}[proposition]{Example}

\def\P{\mathcal{P}}

\def\R{\mathbb{R}}

\def\Q{\mathcal{Q}}

\def\J{\mathcal{J}}

\def\D{\mathcal{D}}

\def\argmin{\mathop{\rm arg\,min}}

\title{\bf{A Lagrangian Dual Based Approach to Sparse Linear
Programming}}

\author{Chen Zhao\thanks{Department of Mathematics, Beijing Jiaotong University, Beijing 100044, P. R. China; 14118409@bjtu.edu.cn.}
\hspace{1mm} 
Ziyan Luo\thanks{Corresponding author; State Key Lab of Rail Traffic Control and Safety, Beijing Jiaotong University, Beijing 100044, P. R. China; starkeynature@hotmail.com.}\hspace{1mm} 
Weiyue Li \thanks{Department of Mathematics, Beijing Jiaotong University, Beijing 100044, P. R. China; 13121542@bjtu.edu.cn.}  \hspace{1mm} 
Houduo Qi \thanks{School of Mathematics, University of Southampton, Southampton SO17 1BJ, UK; hdqi@soton.ac.uk.}\hspace{1mm} 
Naihua Xiu \thanks{Department of Mathematics, Beijing Jiaotong University, Beijing 100044, P. R. China; nhxiu@bjtu.edu.cn.}
}

\begin{document}
\date{April 6, 2018}
\maketitle

\begin{abstract}
A sparse linear programming (SLP) problem is a linear programming problem equipped with a sparsity (or cardinality) constraint, which is nonconvex and discontinuous theoretically and generally NP-hard computationally due to the combinatorial property involved. By rewriting the sparsity constraint into a disjunctive form, we present an explicit formula of its Lagrangian dual in terms of an unconstrained piecewise-linear convex programming problem which admits a strong duality. A semi-proximal alternating direction method of multipliers (sPADMM) is then proposed to solve this dual problem by taking advantage of the efficient computation of the proximal mapping of the vector Ky-Fan norm function. Based on the optimal solution of the dual problem, we design a dual-primal algorithm for pursuing a global solution of the original SLP problem. Numerical results illustrate that our proposed algorithm is promising especially for large-scale problems.
\end{abstract}


\section{Introduction}

A sparse linear programming (SLP) problem is to minimize a linear function subject to a set of linear equality, box and sparsity constraints. 
The involved sparsity constraint
has attracted a great deal of attention in a wide range of applications including compressed sensing \cite{donoho2006compressed}, variable selection \cite{Liu2009Large}, sparse portfolio \cite{Li2006OPTIMAL}, etc. 
Equipped with the sparsity constraint, 
SLP has been employed to reformulate linear compressed sensing \cite{donoho2006compressed}, transportation problems \cite{chen2014sparsity}, Sudoku games \cite{babu2010linear}, and so on. It also resembles the problems of finding the sparsest solution of a linear programming problem \cite{friedlander2007exact,donoho2005sparse} or a linear complementarity problem \cite{shang2014minimal,Chen2016Sparse}, but not exactly so. In contrast to latter cases, SLP is more applicable to the circumstances when people demand a vector with a certain degree of sparsity rather than the sparsest one in real life.

Analogous to most of the zero-norm involved optimization problems, SLP is generally NP-hard due to the intrinsic combinatorial property possessed by the sparsity constraint. 
There are various existing approaches for handling zero-norm related optimization problems, among which two categories can be roughly classified into. The first category consists of approaches with convex and nicely-behaved non-convex relaxation strategies \cite{friedlander2007exact,d2008first}, in which exact relaxation conditions (most of them are not easy to verify in real world applications, such as the restricted isometry property \cite{candes2008restricted}, the null space property \cite{Cohen2009COMPRESSED}, the range space property \cite{Zhao2013RSP}, etc.) are imposed for achieving a sparse solution. 
Approaches in the second category treat the zero-norm directly, including the typical hard-thresholding type algorithms. Most of those approaches are designed for obtaining the sparsest solution and therefore are not suitable for generic SLP.

Existing algorithms for handling sparsity constrained optimization problems in a direct manner are mainly focused on the problems with symmetric sets as feasible regions which admit efficient computation for the corresponding projection operators \cite{beck2015minimization,xu2016efficient,Panlili2017}. 
However, such projected gradient type methods can not be extended to a general SLP problem, since the corresponding feasible region lacks the symmetry and the underlying projection is generally difficult to compute. To our best knowledge, there lacks efficient algorithms for generic SLP, and how to design effect algorithm for solving SLP is still challenging.
This inspires us to develop an effective algorithm for solving such a special type of sparsity constrained problems. 

In this paper, we try to investigate an approach for generic SLP by taking advantage of Lagrangian duality theory which is one of the basic tools for the development of efficient algorithms for mathematical programming \cite{diewert1974applications}. We first give the expression of the Lagrangian dual problem of SLP. 
Then the strong duality for SLP is proved to be true without any assumption. Further, by introducing an auxiliary variable, we reformulate the dual problem as a 2-block separable convex optimization and propose an efficient semi-proximal alternating direction method of multipliers method (sPADMM) to solve it. By developing an efficient searching method for the subproblem involved based on the characterization of the subgradient of the vector Ky-Fan norm function, our sPADMM turns out to be very efficient. Finally, we establish a dual-primal approach for solving SLP. Under some suitable assumptions, we can obtain a global optimal solution of SLP from the solution of dual problem directly.

This paper is organized as follows. The model of SLP and the explicit formula of its Lagrangian dual are stated in Section 2. The strong duality theorem is then established in Section 3, and an sPADMM is proposed for solving the dual problem in Section 4. Based on the sPADMM, a dual-primal algorithm is designed for solving SLP in Section 5. Numerical experiments are conducted in Section 6. Conclusions are made in Section 7.


%

\section{The SLP and Its Lagrangian Dual}
A sparse linear programming (SLP) problem takes the form of
\begin{equation*}
(\P)~~~~\min ~c^Tx,~~\text{s.t.}~Ax=b,~0  \leq  x  \leq l,~\|x\|_0  \leq  r,
\end{equation*}
where $A\in\R^{m\times n}$, $b\in \R^m~c\in\R^n$, $l$ is a positive vector in $\R^n$, $0 < r\leq n$ is a positive integer and $\|x\|_0$ is the $l_0$ norm of $x$, which counts the number of nonzero components in $x$. Assume throughout the paper that the feasible region of $(\P)$ is nonempty, which leads to the solvability of ($\P$) due to the box constraint involved. Note that
$(\P)$ is actually a disjunctive programming since the feasible region can be rewritten as
$$ \bigcup\limits_{I\in \J(r)}\{x\in \R^n_I\mid Ax=b, 0\leq x\leq l\}:= \bigcup\limits_{I\in \J(r)}F(l;I),$$
where 
$\J(r)= \{J\subseteq \{1,\ldots,n\}\mid|J|=r\}$, $\R^n_I$ is the subspace of $\R^n$ spanned by $\{e_i\mid i\in I\}$ with $e_i\in \R^n$ being the $i$-th column of the identity matrix. Thus, $(\P)$ can be cast as
\begin{equation}\label{integer_pro}
  \min\limits_{I\in\J(r)}~\min\limits_{F(l;I)}c^Tx.
\end{equation}
For any given $I\in\J(r)$, the corresponding inner minimization problem
$$(\P_I)~~~~\min\limits_{F(l;I)}c^Tx$$
is a linear program and hence it admits the strong duality. For simplicity, the following notations will be used throughout the paper. Denote the set of all nonnegative (positive) vectors in $\R^n$ by $\R^n_+$ ($\R^n_{++}$). For any given closed set $\Omega\subseteq \R^n$ and any $x\in \R^n$, define $\Pi_\Omega(x):=\argmin\limits_{y\in \Omega} \|x-y\|^2$. The $l_1$ norm, the $l_2$ norm and the projection onto the nonnegative orthant $\R^n_+$ of $x$ are defined as  $\|x\|_1:=\sum\limits_{i=1}^n|x_i|$, $\|x\|_2:=(\sum\limits_{i=1}^n(x_i)^2)^{1/2}$ and $[x]_+:=\max \{x,0\}$, respectively. For any given positive integer $n$, denote $[n]:=\{1,\ldots, n\}$. For any index set $I\subseteq [n]$, $x_I$ represents the subvector of $x$ indexed by index set $I$. Let $x^{\downarrow}$ be the vector of entries of $x$ being arranged in the non-increasing order $x^{\downarrow}_1\geq \dots \geq x^{\downarrow}_n$. The symbol $\circ$ stands for the Hadamard product, i.e., for vectors $u,v\in\R^n$, $(u\circ v)_i=u_iv_i$, $i\in [n]$. For any given optimization problem labeled by (P), we denote its optimal value by Val(P) and its optimal solution set by $S^*_P$.

Given $l\in \R^n_{++}$ and any positive integer $r\leq n$, denote
$S(r):=\{x\in\R^n|\|x\|_0\leq r\}$, $C(l;r):=\{x\in\R^n|0\leq
 x  \leq l,x\in S(r)\}$.
The Lagrangian function of ($\P$) is defined by
$$ L(x,y)=c^Tx-y^T(Ax-b),~~\forall x\in C(l;r), ~\forall y\in \R^m.$$
Then the Lagrangian dual problem takes the form
\begin{equation}\label{dualproblem}
 ~~~ \max\limits_{y\in\R^m}\left\{ \theta(y):=\min\limits_{x\in C(l,r)}L(x,y)\right\}.
\end{equation}
The remainder of this section is devoted to the explicit form of $ \theta(y)$ as defined in \eqref{dualproblem}. Before establishing the desired explicit formula, several essential lemmas are stated for preparation.

\begin{lemma}\label{lemma_opval_I}
Given $I \subseteq [n]$, $p\in\R^n$ and $l\in\R^n_{++}$, we have
\begin{equation}\label{lp}
-\|\Pi_{\R^{n}_I\cap \R^{n}_+}(-l\circ p)\|_1=\min\left\{p^T x: 0\leq x  \leq
 l, x\in \R^n_I\right\}.
\end{equation}
\end{lemma}
\noindent{\bf Proof.}
Set $x^*\in\R^n$ by $x^*_i=l_i$, if $i\in I$ and $p_i<0$, and $x^*_i=0$ otherwise. Then
\begin{equation}\label{ell}
  p^Tx^* = \sum\limits_{i\in I,~p_i<0}l_i p_i
  = -\sum\limits_{i\in I} [-l_i p_i]_+
  = -\|\Pi_{\R^n_I\cap \R^n_+}(-l \circ p)\|_1.
\end{equation}
Note that for any $x\in \R^n_I$ satisfying $0  \leq
 x  \leq
 l$,
$$ p^Tx =  \sum\limits_{i\in I}x_i p_i
  = \sum\limits_{i\in I}x_i ([p_i]_+-[-p_i]_+)
 \geq
  -\sum\limits_{i\in I}x_i [-p_i]_+
\geq
  -\sum\limits_{i\in I}l_i [-p_i]_+
  = p^Tx^*.$$
This completes the proof. \qed


\begin{lemma}\label{dual_I}
The dual problem of ($\P_I$) is
\begin{equation}\label{problem_d_e}
  (\D_I)~~~\max\limits_{y\in\R^m} b^Ty-\|\Pi_{\R^n_I\cap \R^n_+}(l\circ(A^Ty-c))\|_1.
\end{equation}
\end{lemma}
\noindent{\bf Proof.} It follows readily from Lemma \ref{lemma_opval_I} by taking $p=c-A^Ty$.\qed  


\begin{lemma}\label{lemma_spar_pro}
Let $z\in\R^n$ 
and $\pi\in\Pi_{S(r)\cap \R^n_+}(z)$. Then $\pi_{t_i}=[z_{t_i}]_+$ for all $i\in [r]$, and $\pi_{t_i}=0$ otherwise,
where $\{t_1,t_2,\ldots,t_n\}$ satisfies
$z_{t_1}\geq
 z_{t_2}\geq
\ldots\geq
 z_{t_n}$. Furthermore, $\|\pi\|_1=\sum\limits_{i=1}^r [z_{t_i}]_+=\sum\limits_{i=1}^r \left([z]_+\right)^{\downarrow}_i$, for any $\pi\in \Pi_{S(r)\cap \R^n_+}(z).$
\end{lemma}

As indicated in Lemma \ref{lemma_spar_pro}, 
all projection vectors of $z$ onto $S(r)\cap \R^n_+$ share the same $\ell_1$ norm, which allows us to simply use the notation $\|\Pi_{S(r)\cap \R^n_+}(z)\|_1$ for the involved identical $\ell_1$ norm without ambiguity. Combining with Lemma \ref{lemma_opval_I}, this notation can further help characterize the optimal value for a relevant nonconvex optimization problem of ($\P$) as the following lemma states.

\begin{lemma}\label{lemma_opval_S}
Given $p\in\R^n$, $l\in\R^n_{++}$ and an integer $r$ ($0< r\leq
 n$), we have
\begin{equation}\label{pS}
-\|\Pi_{S(r)\cap \R^{n}_+}(-l \circ p)\|_1=\min\left\{p^T x \mid x\in C(l;r)\right\}.
\end{equation}
\end{lemma}
\noindent{\bf Proof.} The desired assertion can be derived by
\begin{eqnarray}\label{eqv}
  \min\limits_{x\in C(l;r)} ~  p^T x &=&\min\limits_{I\in\J(r)} \min\limits_{x\in \R^n_I, 0  \leq
 x  \leq
 l} ~  p^T x  = \min\limits_{I\in\J(r)} -\|\Pi_{\R^n_I\cap \R^n_+}(-l \circ p)\|_1 \nonumber\\
  ~ &=& -\max\limits_{I\in\J(r)}\sum\limits_{i\in I} [-l_i p_i]_+ = -\sum\limits_{i=1}^r \left([-l\circ p]_+\right)^\downarrow_i \nonumber\\
  ~ &=& -\|\Pi_{S(r)\cap \R^{n}_+}(-l \circ p)\|_1,\nonumber
\end{eqnarray} where the first equality follows from the observation $S(r)=\bigcup\limits_{I\in\J(r)} \R^n_I$, the second and the last equalities from Lemmas \ref{lemma_opval_I} and \ref{lemma_spar_pro}, respectively,  and the fourth equality from \eqref{ell}. \qed

Now we are in a position to explicitly formulate the dual problem \eqref{dualproblem}.

\begin{theorem}
The Lagrangian dual problem \eqref{dualproblem} can be explicitly formulated as the following convex unconstrained convex program
$$(\D)~~~~~ \max\limits_{y\in\R^m} b^Ty- \|\Pi_{S(r)\cap \R^n_+}(l\circ (A^Ty-c))\|_1.$$
\end{theorem}
\noindent{\bf Proof.} The formula of $(\D)$ follows readily from the definition of $\theta(y)$ in (\ref{dualproblem}) and Lemma \ref{lemma_opval_S},%
and the convexity follows directly from \cite[Theorem 5.5]{Rockafellar1970Convex}. \qed

%
%
%

\section{Strong Duality}
This section is devoted to the strong duality between $(\P)$ and $(\D)$.

\begin{theorem}\label{strongdual}
\emph{(Strong Duality)} If  ($\P$) (or ($\D$)) is solvable, then so is ($\D$) (or ($\P$)) and the duality gap is zero.
\end{theorem}
\noindent{\bf Proof.} ``$S^*_{\P}\neq \emptyset \Rightarrow S^*_{\D}\neq \emptyset$": Assume that $x^*$ is an optimal solution of ($\P$). The weak duality indicates that
$Val(\D)\leq c^Tx^*$. 
Now we claim that this finite
supremum is attainable. Assume on the contrary that $\theta(y)$ can not reach $Val(\D)$ for any finite $y\in \R^m$. Then there exist  $y_0\in \R^m$ and $d\in\R^m\backslash \{0\}$ 
such that
\begin{equation}\label{limit}
  \lim\limits_{\lambda\rightarrow +\infty}\theta(y_0+\lambda d)=Val(\D).
\end{equation}
Denote $I_+=\{i\in [n] \mid(l \circ A^Td)_i>0\}$, $I_0=\{i\in [n] \mid(l \circ A^Td)_i=0\}$, and $I_-=\{i\in [n] \mid(l \circ A^Td)_i<0\}$.
There always exists a permutation $\{t_1,\ldots, t_n\}$ of [n] and some sufficiently large scalar $\lambda_0>0$ such that
$$\left\{
  \begin{array}{ll}
    (l \circ A^Td)_{t_1}\geq
\cdots\geq
 (l \circ A^Td)_{t_n}, & \hbox{ } \\
    (l\circ (A^T(y_0+\lambda_0d)-c))_{t_1}\geq
\cdots\geq
 (l\circ (A^T(y_0+\lambda_0d)-c))_{t_n}, & \hbox{ } \\
    (l\circ (A^T(y_0+\lambda_0d)-c))_i>0,~\forall i\in I_+, & \hbox{ } \\
    (l\circ (A^T(y_0+\lambda_0d)-c))_k<0,~\forall k\in I_-. & \hbox{ }
  \end{array}
\right.$$
By direct calculation, we can further get that for any $\alpha>0$,
\begin{eqnarray}\label{ordering}
&&\|\Pi_{S(r)\cap \R^n_+}(l\circ(A^T(y_0+\lambda_0d)-c)+\alpha l\circ A^Td)\|_1\nonumber\\
&=&\|\Pi_{S(r)\cap \R^n_+}(l\circ (A^T(y_0+\lambda_0d)-c))\|_1+\alpha \|\Pi_{S(r)\cap \R^n_+}l\circ A^Td)\|_1.
\end{eqnarray}
Thus, for any $\lambda\geq
\lambda_0$,
\begin{eqnarray*}
\theta(y_0+\lambda d) &=&b^T(y_0+\lambda d)-\|\Pi_{S(r)\cap \R^n_+}(l\circ (A^T(y_0+\lambda d)-c))\|_1\\
  ~ &=&  b^T(y_0+\lambda_0 d)+(\lambda-\lambda_0)b^Td-\|\Pi_{S(r)\cap \R^n_+}(\lambda-\lambda_0)l\circ A^Td)\|_1\\
  ~ &~&  -\|\Pi_{S(r)\cap \R^n_+}(l\circ (A^T(y_0+\lambda_0d)-c))\|_1 \\
  ~ &=& \theta(y_0+\lambda_0 d)+(\lambda-\lambda_0)(b^Td-\alpha\|\Pi_{S(r)\cap \R^n_+}l\circ A^Td)\|_1),
\end{eqnarray*}
where the second equality follows from \eqref{ordering}. Combining with (\ref{limit}), we have $\theta(y_0+\lambda_0 d)=Val(\D)$, which is a contradiction.
Thus, ($\D$) is solvable.

\noindent``$S^*_{\D}\neq \emptyset \Rightarrow S^*_{\P}\neq \emptyset$": Let
$y^*$ be an optimal solution of $(\D)$. Then
there exists some index set $\hat{I}\in \J(r)$ such that
$y^*\in \arg\max\limits_{y\in\R^m} b^Ty- \|\Pi_{\R^n_{\hat{I}}\cap \R^n_+}(l\circ (A^Ty-c))\|_1.$ By invoking Lemma \ref{dual_I}, together with the strong duality theorem of linear programming, 
we can find some $\hat{x}\in F(l;{\hat{I}})$ such that
\begin{equation}\label{uu} Val(\P_{\hat{I}})=c^T\hat{x}= b^Ty^*- \|\Pi_{\R^n_{\hat{I}}\cap \R^n_+}(l\circ (A^Ty^*-c))\|_1=Val(\D).\end{equation}
Such a vector $\hat{x}$ is certainly a feasible solution of $(\P)$ and hence $(\P)$ is solvable due to the compactness of the feasible region of $(\P)$.

The remainder is to show the zero duality gap between $(\P)$ and $(\D)$. Assume that $x^*$ is an optimal solution of $(\P)$ and $y^*$ is an optimal solution of $(\D)$. It follows from \eqref{integer_pro} that there exists an index set $I^*\in\J(r)$ such that $x^*$ is an optimal solution of ($\P_{I^*}$) and for any $I\in\J(r)$, we have \begin{equation}\label{vv}
  Val(\P_I)\geq
 Val(\P_{I^*})
\end{equation}
Suppose that the duality gap is positive. Then $ Val(\P_{I^*})=c^Tx^* > Val(\D)= Val(\P_{\hat{I}})$.
where the last equality follows readily from \eqref{uu}. This is a contradiction to \eqref{vv}. This completes the proof. \qed

\section{Solving the Dual Problem via sPADMM}
As established in Section 2, the dual problem $(\D)$ is a non-smooth convex problem with a piecewise-linear objective function. By introducing an auxiliary variable vector $z$, the dual problem can be rewritten as
$$(\hat{\D}) ~~~~\min ~  f(y)+h(z):=-b^Ty+\|\Pi_{S(r)\cap \R^n_+}(l\circ z)\|_1~~~~ \text{s.t.}~~ A^Ty-z=c.$$
Since problem ($\hat{\D})$ is a 2-block separable convex optimization problem, the semi-proximal alternating direction method of multipliers (sPADMM) as introduced by Sun et al. in  \cite{sun2015convergent} can be applied to solve ($\hat{\D})$.

For any given penalty parameter $\sigma>0$ and $(y,z,w)\in \R^m\times \R^n\times \R^n$, the augmented Lagrangian function of ($\hat{\D})$ is defined by
\begin{equation*}
  L_\sigma(y,z;w)=-b^Ty+\|\Pi_{S(r)\cap \R^n_+}(l\circ z)\|_1-w^T(A^Ty-z-c)+\frac{\sigma}{2}\|A^Ty-z-c\|^2.
\end{equation*}
By choosing an initial point $(y^0,z^0,w^0)$, the sPADMM takes the following iteration scheme for $k=0,1,\ldots$,
\begin{equation}\label{scheme}
\left\{
  \begin{array}{ll}
    y^{k+1}=\argmin  L_\sigma(y,z^k;w^k)+\frac{\sigma}{2}\|y-y^k\|^2_{\P}, & \hbox{ } \\
   z^{k+1}= \argmin  L_\sigma(y^{k+1},z;w^k)+\frac{\sigma}{2}\|z-z^k\|^2_{\Q}, & \hbox{  } \\
    w^{k+1}=w^k-\tau\sigma(A^T y^{k+1}-z^{k+1}-c), & \hbox{ }
  \end{array}
\right.
\end{equation}
where $\P$ and $\Q$ are two symmetric positive semidefinite matrices of appropriate sizes, and $\|x\|_\P:=\sqrt{x^T\P x}$ is the $\P$-weighted norm. The choices of $\P$ and $\Q$ should in principle as small as possible and meanwhile to simplify the computation of the corresponding subproblems in the aforementioned scheme. Details for solving ($\hat{\D}$) by \eqref{scheme} will be elaborated as follows.

\subsection{The Update of $y$}
Given the current iteration point $(y^k,z^k,w^k)\in \R^m\times \R^n\times \R^n$, the update of the $y$-part is obtained by
\begin{eqnarray}\label{newy}
 y^{k+1} &=&  \argmin  L_\sigma(y,z^k;w^k)+\frac{\sigma}{2}\|y-y^k\|^2_{\P}, \nonumber \\
   &=& \left(AA^T+\P\right)^{-1}\left(A(w^k+\sigma z^k+\sigma c)+b+\sigma \P y^k\right).
\end{eqnarray}
Evidently, if $AA^T$ is nonsingular and admits a fast way to get its inverse, then we can simply choose $\P=0$. Otherwise, we can take $\P=\lambda_{\max}(AA^T)-AA^T$ to ease the computation of $y^{k+1}$ via $y^{k+1}=\frac{1}{\sigma\lambda_{\max}(AA^T)}(A(w^k+\sigma z^k+\sigma c)+b+\sigma\P y^k)$, where $\lambda_{\max}(AA^T)$ stands for the largest eigenvalue of $AA^T$.

\subsection{The Update of $z$}
Given $(y^{k+1},z^k,w^k)\in \R^m\times \R^n\times \R^n$, the update of $z$ is obtained by
\begin{eqnarray}\label{newz0}
z^{k+1} &=&  \argmin L_\sigma(y^{k+1},z;w^k)+\frac{\sigma}{2}\|z-z^k\|^2_{\Q}, \nonumber \\
   &=&  L^{-1}\argmin  \|\Pi_{S(r)\cap \R^n_+}(z)\|_1+\frac{\sigma}{2}z^TMz-\sigma\left\langle z, \tilde{w}\right\rangle
\end{eqnarray}
where $L:=Diag(l)$, $M:=L^{-1}(I+Q)L^{-1}$ and $\tilde{w}:=L^{-1}\left(A^Ty^{k+1}-c+Qz^k-\frac{w^k}{\sigma}\right)$. It is easy to find that if $l=l_0 e$ for some positive scalar $l_0$, then we can simply choose $\Q=O$ and update $z$ by
\begin{equation}\label{newz1}
z^{k+1}=\argmin \frac{1}{2}\|z-\bar{w}\|^2 + \bar{\lambda} \|\Pi_{S(r)\cap \R^n_+}(z)\|_1,
\end{equation} with $\bar{w}=(A^T y^{k+1}-c)- \frac{w^k}{\sigma}$ and $\bar{\lambda}=\frac{l_0}{\sigma}$.
Otherwise, we can choose $\Q=\frac{L^2}{l^2_{\min}}-I$ and update $z$ by
\begin{equation}\label{newz2}
z^{k+1}= L^{-1}\left(\argmin \frac{1}{2}\|z-\tilde{w}\|^2 + \tilde{\lambda} \|\Pi_{S(r)\cap \R^n_+}(z)\|_1\right)
\end{equation}  with $\tilde{w}$ as defined in \eqref{newz0} and $\tilde{\lambda}:=\frac{l^2_{\min}}{\sigma}$.
Note that in both cases of \eqref{newz1} and \eqref{newz2}, the essential task for computing $z^{k+1}$ is to handle the following minimization problem
\begin{equation}\label{reg1}
  \min~~\frac{1}{2}\|z-w\|^2+\lambda\|\Pi_{S(r)\cap\R^n_+}(z)\|_1
\end{equation}
for some given $w\in \R^n$, $\lambda>0$ and $0<r\leq
 n$. Thus, the reminder of the subsection is devoted to proposing efficient algorithms for solving \eqref{reg1}.

For any given $w\in \R^n$, select a permutation $\pi$ of $[n]$ such that $w^\downarrow=w_{\pi}$.
 We can verify that \eqref{reg1} is equivalent to \begin{equation}\label{reg11}
  \min~~\frac{1}{2}\|z-w^\downarrow\|^2+\lambda\|\Pi_{S(r)\cap\R^n_+}(z)\|_1
\end{equation}
in the sense that $z^*$ solves problem (\ref{reg11}) if and only if $z^*_{\pi^{-1}}$ solves problem (\ref{reg1}).
Divide $w^\downarrow$ into two parts,
\begin{equation}\label{w+-}
  w^\downarrow=({w^+}^T,{w^-}^T )^T
\end{equation}
with $w^+\in \R^{n_1}_+$ and $w^-\in (-\R^{n_2}_{++})$ with $n_1+n_2=n$. We now consider the following low-dimensional subproblem
\begin{equation}\label{reg4}
  \min\limits_{z\in\R^{n_1}}~~g(z):=\frac{1}{2}\|z-w^+\|^2+\lambda\|z\|_{(r)},
\end{equation}
where $\|z\|_{(r)}=\sum\limits_{i=1}^r|z|^\downarrow_i$, $w^+=\left([w]_+\right)^\downarrow$ are the same as defined in \eqref{w+-}, $\lambda>0$ and $0<r\leq
 n_1$. Apparently, the unique solution of \eqref{reg4} is the value of the proximal mapping of the vector Ky-Fan $r$-norm function $\|\cdot\|_{(r)}$ at $w^+$, 
say $\bar{z}$, which satisfies the following necessary and sufficient optimality condition
\begin{equation}\label{sear2}
0\in \partial g(\bar{z})=\bar{z}-{w^+}+\lambda \partial \|\bar{z}\|_{(r)}.
\end{equation}
Two useful lemmas are reviewed for better characterizing $\bar{z}$.

\begin{lemma}[Theorems 368,369,~\cite{hardy1952inequalities}]\label{lemma_ineq}
For any $x$, $y \in\R^n$, $\langle x,y \rangle\leq\langle x^\downarrow,y^\downarrow \rangle$,
where the inequality holds if and only if there exists a permutation $\pi$ of $[n]$ such that $x_\pi=x^\downarrow$ and $y_\pi=y^\downarrow$.
\end{lemma}

\begin{lemma}[\cite{overton1993optimality,watson1993matrix}]\label{sear3}
Assume that $z\in \R^n_+$ satisfies 
\begin{equation*}
  z_1\geq
\ldots\geq
 z_{r_0}>z_{r_0+1}=\ldots=z_r=\ldots=z_{r_1}>z_{r_1+1}\geq
\ldots\geq
 z_n\geq
0,
\end{equation*}
where $r_0$ and $r_1$ are integers such that $0\leq
 r_0<r \leq
 r_1\leq
 n$. If $z_r>0$, then
\begin{equation*}
  \partial \|z\|_{(r)}=\left\{\mu \left|
  \begin{split}
  \mu_i=1,i=1,\ldots,r_0,~\mu_i=0,i=r_1,\ldots,n\\
  0\leq
 \mu_i \leq
 1, \sum\limits^{r_1}_{j=r_0+1}\mu_j=r-r_0,i=r_0+1,\ldots,r_1
  \end{split}\right.\right\}.
\end{equation*}
Otherwise,
\begin{equation*}
  \partial \|z\|_{(r)}=\left\{\mu  \left|
  \mu_i=1,i=1,\ldots,r_0,~
  0\leq
 \mu_i \leq
 1, \sum\limits^{n}_{j=r_0+1}\mu_j\leq
 r-r_0,i=r_0+1,\ldots,n\right.\right\}.
\end{equation*}

\end{lemma}

\begin{proposition}\label{pr1}
$\bar{z}$ solves Problem \eqref{reg4} with $\bar{z}_r=0$ if and only if

\begin{equation}\label{kkt2}
\left\{
              \begin{array}{ll}
                \bar{z}_i=w^+_i-\lambda, &i=1,\ldots,r_0 \\
                \bar{z}_i=0, &i=r_0+1,\ldots,n
              \end{array}
            \right.
\end{equation}
 where the index $r_0$ satisfies $w^+_{r_0}>\lambda \geq
 w^+_{r_0+1}$ and $r-r_0\geq
 \sum\limits^{n}_{j=r_0+1}w^+_j/\lambda$.

\end{proposition}

\noindent{\bf Proof.} Note that $\bar z =|\bar z|^\downarrow$. Thus, $\bar{z}_r=0$ if and only if there exists an index $r_0\in[r-1]$satisfying
\begin{equation}\label{indexr0}
  \bar z_1\geq
\ldots\geq
 \bar z_{r_0}>\bar z_{r_0+1}=\ldots=\bar z_r=\ldots=\bar z_{r_1}=\ldots \bar z_n = 0.
\end{equation}
According to Lemma \ref{sear3} the optimal condition (\ref{sear2}) is equivalent to
\begin{equation}\label{kkt}
\left\{
  \begin{array}{ll}
    \bar{z}=w^++\lambda \mu,\\
    \bar z_{r_0}>\bar z_{r_0+1}=\ldots=\bar z_r=\ldots=\bar z_{r_1}=\ldots \bar z_n = 0,\\
    \mu_i=1,i=1,\ldots,r_0,~
  0\leq
 \mu_i \leq
 1, \sum\limits^{n}_{j=r_0+1}\mu_j\leq
 r-r_0,i=r_0+1,\ldots,n.\\
  \end{array}
\right.
\end{equation}
By solving (\ref{kkt}), we obtain
\begin{equation}\label{kkt2}
\left\{
              \begin{array}{ll}
                \bar{z}_i=w^+_i-\lambda, &i=1,\ldots,r_0 \\
                \bar{z}_i=0, &i=r_0+1,\ldots,n
              \end{array}
            \right.
\end{equation}
By invoking the structure of $|v|^\downarrow$, we obtain the equivalence between (\ref{kkt}) and
$$w^+_{r_0}>\lambda \geq
 w^+_{r_0+1},~ r-r_0\geq
 \sum\limits^{n}_{j=r_0+1}w^+_j/\lambda.$$
This completes the proof.\qed 

\begin{proposition}\label{pr2}
$\bar{z}$ solves Problem \eqref{reg4} with $\bar{z}_r>0$ if and only if

\begin{equation}\label{kkt4}
\left\{
              \begin{array}{ll}
               \bar{z}_i= w^+_i-\lambda, &i=1,\ldots,r_0 \\
                \bar{z}_i=\theta, &i=r_0+1,\ldots,r_1 \\
                \bar{z}_i=w^+_i, &i=r_1+1,\ldots,n
              \end{array}
            \right.
\end{equation}
where the index $r_0$ satisfies $w^+_{r_0}>\lambda+\theta \geq
 w^+_{r_0+1}$ and $w^+_{r_1}\geq
\theta > w^+_{r_1+1}$ with $\theta= (\sum\limits^{n}_{j=r_0+1}w^+_j-\lambda(r-r_0))/(r_1-r_0)$.

\end{proposition}

\noindent{\bf Proof.} By mimicking the proof of Proposition \ref{pr1}, we can also get the desired results. Thus we omit the proof here.\qed


By employing Propositions \ref{pr1} and \ref{pr2}, we can adopt a searching method proposed by \cite{wu2014moreau} to compute the proximal mapping induced by $\|\cdot\|_{(r)}$ exactly and hence get the desired $\bar{z}$ for Problem \eqref{reg4}.
Define $s\in \R^{n_1+1}$ by $s_0=0$ and $s_j=\sum\limits^j_{i=1} w_i^+$, $j=1,\ldots,n_1$.
Let $v^-$ and $v^+$ be two vectors of length $r+1$ and $n_1-r+2$ such that $v^-_0=+\infty$, $v^-_i=w^+_i$, for any $i\in [r]$, and $v^+_{p+1}=0$, $v_i^+=w^+_i$, for any $i\in [r]$.
Algorithm 1 is proposed to get the optimal solution $\bar{z}$ of Problem (\ref{reg4}).

\begin{algorithm} \label{searching}
\caption{Searching($\lambda, w^+, r$)}

{\bf step 0.} Pre-compute $s,~v^-,~v^+$, set $r_0=r-1$.

{\bf Step 1.} If $w^+_{r_0}>\lambda \geq
 w^+_{r_0+1}$ and $\lambda\geq
 \sum\limits^{n}_{j=r_0+1}w^+_j/(r-r_0)$, compute $\bar{z}$ by \eqref{kkt2}
and go to {\bf Step 3}. Otherwise if $r_0=0$, set $r_0=r-1$, $r_1=r$ and go to {\bf Step 2}; if $r_0>0$, replace $r_0$ by $r_0-1$ and repeat {\bf Step 1}.

{\bf Step 2.} Set flag=0 and  compute $\theta=(\sum\limits^{n}_{j=r_0+1}w^+_j-\lambda(r-r_0))/(r_1-r_0)$. If $r_0=0$ and $r_1=n_1$, set flag=1. Otherwise if $w^+_{r_0}>\lambda+\theta \geq
 w^+_{r_0+1}$ and $
  w^+_{r_1}\geq
\theta > w^+_{r_1+1}$, set flag=1. If flag=1, compute $\bar{z}$ by \eqref{kkt4} and go to {\bf Step 3}. If flag=0, and $r_1<n_1$, replace $r_1$ by $r_1+1$ and repeat {\bf Step 2}; otherwise replace $r_0$ by $r_0-1$, set $r_1=r$ and repeat {\bf Step 2}.

{\bf Step 3.} Output $\bar z$.

\end{algorithm}

%
%

Based on the optimal solution $\bar{z}$ of Problem (\ref{reg4}), we can then get the optimal solution of the subproblem \eqref{reg11} as presented in the following theorem and finally get the optimal solution of the subproblem \eqref{reg1} by permutation.

\begin{theorem}  For any given $w\in \R^n$ with the decomposition \eqref{w+-}, let $\bar{z}$ and $z^*$ be the unique solution of Problem (\ref{reg4}) and Problem (\ref{reg11}), respectively. Then
$z^*=( {\hat{z} }^T,{w^- }^T  )^T$, where $\hat{z_i}=\left\{
  \begin{array}{ll}
    w^+_i, & \mbox{~if~}w^+_i>\lambda \\
    0, & \mbox{~if~} 0\leq
 w^+_i \leq
 \lambda
  \end{array}
\right.$ for each $i\in [n_1]$ if $r\geq
 n_1$, and $\hat{z}=\bar{z}$ otherwise.
\end{theorem}
\noindent{\bf Proof.} By virtue of \eqref{w+-} and (\ref{reg11}), we can verify that
$\left(z^*\right)_{i=n_1+1}^n =w^-$, and $\left(z^*\right)_{i=1}^{n_1} = \arg\min\limits_{z\in\R^{n_1}_+} \frac{1}{2}\|z-w^+\|^2+\lambda\|\Pi_{S(r)}(z)\|_1.$
Evidently, for any $z\in\R^{n_1}_+$, $\|\Pi_{S(r)}(z)\|_1=\|z\|_{(r)}$. Moreover, if $r\geq n_1$, then $\|\Pi_{S(r)}(z)\|_1 = e^Tz$. This further implies the desired assertion.\qed

\subsection{The sPADMM Algorithm}
The framework of the sPADMM for solving ($\hat{\D})$ is stated in Algorithm 2.
\begin{algorithm}[ht]\label{sPADMM1}
\caption{An sPADMM for Solving ($\hat{\D})$ }

Choose $\sigma>0$ and $\tau\in(0,\frac{1+\sqrt{5}}{2})$ as parameters, and $\left(y^0,~z^0,~w^0\right)$ as the initial point. For $k=0,1,\ldots$, do the following iteration unless some stopping criteria are met:

{\bf Step 1.} Compute $y^{k+1}$ by \eqref{newy};

{\bf Step 2.} Compute $z^{k+1}$ by \eqref{newz1} if $l=l_0 e$ for some $l_0>0$, or by \eqref{newz2} otherwise, where the optimal solution $z^*$ of involved subproblem computed in the following way:

{\bf If $n_1\leq r$,} then for each $i\in [n_1]$, set $\bar z_i=w^+_i$ if $w^+_i>\lambda$, and $\bar z_i=0$, if $0\leq w^+_i \leq  \lambda$;
%

{\bf Else} $\bar z=searching(\lambda,w^+,r).$  Compute $w^-$ to get $z^*=\left({ \bar z}^T,{w^-}^T\right)^T$;

{\bf Step 3.} $w^{k+1} = w^k-\tau\sigma(A^T y^{k+1}-z^{k+1}-c)$.

\end{algorithm}
Note that $dom f=\R^m$, $dom h=\R^n$ and $(0,-c)\in \{(y,z) \mid A^T y-z=c\}$ in Problem ($\hat{\D})$. The following convergence result of Algorithm 2 follows readily from \cite[Theorem B.1]{fazel2013hankel}.
%
%
%
\begin{theorem}
 Assume that ($\hat{\D})$ is solvable. Let $\{(y^k,z^k,w^k)\}$ be generated from Algorithm 2. If $\tau\in(0,(1+\sqrt{5})/2)$, then the sequence $\{(y^k,z^k)\}$ converges to an optimal solution of ($\hat{\D})$ and $\{w^k\}$ converges to an optimal solution to the dual problem of ($\hat{\D}$).
\end{theorem}

The dual problem of ($\hat{\D}$) admits an explicit form as follows.

\begin{lemma}\label{d-d}
  The dual problem of ($\hat{\D}$) can be formulated as
  \begin{equation*}
(\hat{\P})~~~~\min~~c^T w~~~\text{s.t.}~~  Aw=b, ~ w^Tl^{-1}\leq 1,~w\geq 0. 
\end{equation*}

\end{lemma}

\noindent {\bf Proof.} The explicit formula of the dual problem of ($\hat{\D}$) follows readily from Lemma \ref{lemma_opval_S}, \cite[Example 11.4]{VA} and the observation that $cl(conv(C(l;r))) = \left\{w\in\R^n\left|w^T l^{-1}\leq 1, w\geq 0\right.\right\}$. \qed 

As indicated in Lemma \ref{d-d}, the feasible region of the dual problem of ($\hat{\D}$) is a subset of $\left\{x\in\R^n_+\mid Ax=b,x\geq l\right\}$. By utilizing the strong duality as shown in Theorem \ref{strongdual}, we can verify that $S_{\P}^* = S_{\hat{\P}}^*\cap S(r)$. 
In this case, we can get a global solution of $(\P)$ via Algorithm 2 once the sparsity constraint is satisfied. Otherwise, a dual-primal approach will be proposed to handle those non-sparse cases, which will be discussed in the next section.



\section{A Dual-Primal Approach}

In this section, a dual-primal algorithm will be designed to get a global solution of the original sparse linear programming problem $(\P)$ based on Algorithm 2. Before proceeding, we introduce the concept of optimal index sets of $(\P)$.
\begin{definition}
An index set $I^*\in \J(r)$ is said to be an optimal index set of $(\P)$ if there exists an optimal solution $x^*$ such that $x^*\in\R^n_{I^*}$.
\end{definition}

The following theorem provides several sufficient conditions for getting an optimal index set of $(\P)$ by utilizing the optimal solution of ($\hat{\D}$).

\begin{theorem}\label{thindex}
Suppose that $(y^*,z^*)$ is the optimal solution of ($\hat{\D}$). Let $\{t_1,\ldots,t_n\}$ be a permutation of $[n]$ such that 
$l_{t_1} z^*_{t_1}\geq l_{t_2} z^*_{t_2}\geq\ldots\geq l_{t_n} z^*_{t_n}$.
Then,
\begin{itemize}
  \item[(i)] $0$ is an optimal solution of $(\P)$ if $z^*<0$; 
  \item[(ii)] $\{t_1,\ldots,t_r\}$ is an optimal index set of $(\P)$ if one of the following holds
  \begin{itemize}
  \item[(a)] $\|[z^*]_+\|_0=r$;
  \item[(b)] $\|[z^*]_+\|_0>r$ and $l_{t_r}z^*_{t_r}>l_{t_{r+1}}z^*_{t_{r+1}}$;
  \item[(c)] $0<\|[z^*]_+\|_0<r$ and $z^*_{t_{r+1}}<0$.
\end{itemize}
\end{itemize}
\end{theorem}

\noindent{\bf Proof.} 
Note that $y^*$ is also the optimal solution to $(\D)$. Thus,
there exists an index set $I^*\in \J(r)$ such that
 $$b^Ty^*-\|\Pi_{\R^n_{I^*}\cap \R^n_+}(l\circ(A^Ty^*-c))\|_1=b^Ty^*-\|\Pi_{S(r)\cap \R^n_+}(l\circ(A^Ty^*-c))\|_1.$$
It then follows from Lemma \ref{lemma_spar_pro} that for any index set $I\in\J(r)$,
 \begin{eqnarray}
   \|\Pi_{\R^n_{I^*}\cap \R^n_+}(l\circ(A^Ty^*-c))\|_1 &\geq& \|\Pi_{\R^n_{I}\cap \R^n_+}(l\circ(A^Ty^*-c))\|_1.\label{II}
 \end{eqnarray}
Let $x^*$ be an optimal solution of ($\P_{I^*}$).
It follows from the proof of Theorem \ref{strongdual} that $x^*$ is exactly an optimal solution of $(\P)$. By invoking \eqref{II}, we can get the assertions in (a) and (b) immediately. From the first-order optimality for ($\P_{I^*}$), we know that $(x^*,y^*)$ will satisfy 
$$\left\{
  \begin{array}{ll}
    (A^Ty^*-c)_{I^*}=-u+v,\\
    u^Tx^*_{I^*}=0,~x^*_{I^*}\geq0,~u\geq0,\\
    v^T(l_{I^*}-x^*_{I^*})=0, ~l-x^*_{I^*}\geq0,~v\geq0.
  \end{array}
\right.$$
Note that $z^*=A^Ty^*-c$. By direct calculation, we can verify that for any $i\in I^*$, if $z^*_i<0$, then $x^*_i=0$. Thus, it is not hard to verify the results in (i) and (c) are valid. \qed

Theorem \ref{thindex} allows us to design a dual-primal algorithm for pursuing a global solution of $(\P)$ by solving $(\hat{\D})$ and a corresponding linear programming problem in an $r$-dimensional subspace indexed by the optimal index set of $(\P)$, supposing that one of the conditions in (a) and (b) of Theorem \ref{thindex}. The algorithm is stated in Algorithm 3.

\begin{algorithm}
\caption{A Dual-Primal Algorithm for ($\P$)}
{\bf Step 1.}(Dual) Solve ($\hat{\D}$) by Algorithm 2 to obtain $(y^*,w^*)$.

{\bf IF} $\|w^*\|_0\leq r$, $x^*=w^*$;
{\bf Else} go to Step 2;

{\bf Step 2.}(Primal) Compute $z^*=A^Ty^*-c$. Sort $z^*$ into a non-increasing order with a permutation $\pi$. Set $I^*=[\pi_1,\ldots,\pi_r]^T$. Solve Problem ($\P_{I^*}$) to obtain $x^*$.

\end{algorithm}

%
%
%
%
%

It is worth pointing out that we might fail to get the optimal index set $I^*$ once none of the conditions in Theorem \ref{thindex} is satisfied for the optimal solution $(y^*, z^*)$ of $(\hat{\D})$, as the following example illustrates.

\begin{example}\label{exindex} Set
 $A= \left(
    \begin{array}{cccc}
      1 & -1 & 0 & 0 \\
      0 &  0 & 1 & -1 \\
    \end{array}
  \right)$, $b=(0,0)^T,~c=(-1,-1,-1,-1)^T$, $l=(1,1,1,1)^T$ and $r=2$ in problem $(\P)$. One can find that the optimal solution set of ($\P$) is $\{(1,1,0,0)^T, (0,0,1,1)^T\}$ with the
optimal value $-2$. By simple calculation, we know that the optimal solution of its dual problem is $(y^*,z^*)$ with $y^*=(0,0)^T$ and  $z^*=(1,1,1,1)^T$. As can be seen that none of the conditions on $(y^*,z^*)$ presented in Theorem \ref{thindex} hold, and hence no optimal index set of $(\P)$ can be identified.
\end{example}

In the case as stated in Example \ref{exindex}, Algorithm 3 fails to work. This is reasonable due to the NP-hard complexity of $(\P)$. How to handle such a case remains an open question.



\section{Numercial Experiments}

In this section, we will evaluate the performance of our proposed dual-primal algorithm (DPA) for solving sparse linear programming problems. 
As a main part of DPA, the quality of the solution generated by sPADMM will greatly affect the accuracy of the solution generated by DPA. In our numerical experiments, we will measure the accuracy of an approximate optimal solution $\left(\tilde{y},\tilde{z}, \tilde{w}\right)$ of $(\hat{\D})$ and $(\hat{\P})$ by using the following relative infeasibility and the relative duality gap:
$$\zeta:= \frac{A^T\tilde{y}-\tilde{z}-c}{\|c\|}, ~~\eta := \frac{c^T \tilde{w}-b^T\tilde{y}}{\max\{|c^T\tilde{w}|,b^T\tilde{y}\}}.$$
For a given tolerance $\epsilon>0$, we will stop the algorithm when $\zeta<\epsilon$ and $\eta<\epsilon$. For all the tests here, we set $\epsilon = 1$e-8.
The algorithm will also be stopped when it reaches the maximum number of iterations, say Maxiter. Here we set Maxiter = 5e3. We simply choose the dual step-length $\tau=1.618$. 
All the computational results are obtained from 
a desktop computer running on 64-bit Windows Operating System having 4 cores with Intel(R) Core(TM) i5-5257U CPU at 2.70GHz and 16 GB memory.

\subsection{Randomly Generated Examples}
We first test some randomly generated examples whose $I^*$ can be obtained through $y^*$ via Theorem \ref{thindex}. We generate the ``true" optimal solution $xopt$ by using the MATLAB commands ``$xopt = zeros(n,1)$; $T = randperm(n)$; $r = ceil(rand*r)$; $xopt(T(1:r)) = abs(randn(r,1));$"
with different choices of $n$ and $r$. The input data $A\in\R^{m\times n}$, $b\in\R^m$, $l\in\R^n_{++}$ and $c\in\R^n$ are generated by the MATLAB commands ``$A = randn(m,n)$; $b = A*xopt$; $alpha = max(xopt)$; $l = alpha*ones(n,1)$; $c = ones(n,1)$; $c(xopt>0) = 0;$".
Note that Problem $(\P)$ can be cast as a mixed integer programming (MIP) problem
$$ \min\limits_{x\in\R^n_+}\left\{c^Tx\mid Ax=b,x\leq l, z\in \{0,1\}^n, (e-z)+eps\leq x\leq Mz, e^Tz\leq r\right\},$$
where $eps>0$ is sufficiently small (e.g., $eps$=1e-10), and $M>0$ is sufficiently large (e.g., $M$=1e15). Thus, we can use cvx by calling Gurobi to solve this MIP problem.
%
%
%
%
Set $n=1000$, $m=500$ and  $r=25:5:50$. For all those six choices of $r$, we run 100 independent instances to get an average performance of our proposed DPA and the MIP solver (Gurobi). Results are reported in Table 1. The ``Success rate" of DPA in Table 1 stands for the percentage of successful results of 100 instances in the sense that if the relative error of $x$ is smaller than ${10}^{-2}$, i.e., $\frac{\|x-xopt\|}{\|x\|}<{10}^{-2}$,
the result is regarded as a success. It is known from Table 1 that for such a special structured class of SLP instances, both DPA and MIP are robust with the sparsity $r$ in terms of computation time. And by taking the advantage of the special structure of the input data, our DPA is faster than the MIP approach. %
\begin{table}[h]
\caption{Randomly generated examples with $n=1000$ and $m=500$}
\vspace{.1 in}
\begin{center}
  \begin{tabular}{c|ccc|c}
  \hline
    &                        & DPA                           &                       &       MIP   \\
    \textbf{r}              & \textbf{ Iter No. }          & \textbf{ Success rate }& \textbf{CPU time (s)}& \textbf{CPU time (s)}\\
 \hline
 10                    &  84                          & 100\%         &         0.3251  &      7.5010      \\

 25                    &  92                         & 100\%         &         0.3474  &      7.5512           \\

 50                    &  95                         & 100\%        &         0.3538  &      7.1546           \\

100                    &  113                        & 100\%         &        0.4204 &       7.2447             \\
 \hline
\end{tabular}
\end{center}
\end{table}

For relatively large-scale instances, we report the average numerical results of DPA and MIP with $n$ varying from $5000$ to $10000$, $m=0.3n$ and $r=0.05n$ in Table 2 by running each case with 50 independent instances.
\begin{table}[h]
\caption{50 independent examples with $n$ increase from $5000$ to $10000$}
\vspace{.1 in}
\begin{center}
   \begin{tabular}{c|ccc|c}
  \hline
    &                        & DPA                           &                       &       MIP   \\
    \textbf{r}     & \textbf{ Iter No. } & \textbf{ Success rate }& \textbf{CPU time (s)}& \textbf{CPU time (s)}\\
 \hline
 5000              &   170                &  100\%                 & 10.3380               & 104.2140\\
 7000              &   184                &  100\%                 & 27.0782               & 224.4301\\
 9000              &   225                &  100\%                 & 29.7644               & 364.7010\\
 10000             &   279                &  100\%                 & 53.5050               & 467.1680\\
 \hline
\end{tabular}
\end{center}
\end{table}

As one can see in both Tables 1 and 2, our proposed DPA can provide a global solution of each instance and outperforms the MIP approach in the above special setting.

\subsection{Simplex-Constrained SLP}
In \cite{xu2016efficient}, Xu et al. proposed a nonmonotone projected gradient (NPG) method for solving the following simplex-constrained sparse optimization problem
$$\min ~f(x),~~\text{s.t.}~e^Tx=1,~0 \leq x \leq ue,~\|x\|_0 \leq r.$$
where $f:\R^n\rightarrow \R$ is Lipschitz continuously differentiable. If we take $f(x) = c^T x$, it becomes a special case of our model with $A = e^T$, $b = 1$ and $l = ue$. Our algorithm can also solve this model efficiently.
For simplicity, the inputs are generated by using the Matlab commands ``$c = randn(n,1); u = 1$".

By running 100 instances with $n$ varying from 5000 to 10000 and $r=0.05n$, the average performance of DPA is shown in Table 3, in which the ``success rate" states that how many percentage of the tested instances on which DPA can provide an approximate solution within the relative error ${10}^{-2}$ between that generated by NPG. As one can see in Table 3, all instances have been successfully solved within half a second.



\begin{table}[ht]
 \caption{100 independent examples with different choices of $n$}
 \vspace{.1 in}
\begin{center}
   \begin{tabular}{c|ccc}
  \hline
  \textbf{r}     & \textbf{ Iter No. } & \textbf{ Success rate }& \textbf{CPU time (s)} \\
 \hline
 5000            &         249   &   $100\%$ &   0.1163     \\

 7000            &         332   & $100\%$   &  0.1478      \\

 9000            &         454   &  $100\%$  &  0.3987     \\

 10000           &         455   &  $100\%$  &  0.4411     \\

 \hline

\end{tabular}
\end{center}
\end{table}

%
%
%
%
%

For relatively large-scale cases with different sparsity constraints, e.g., $n = 20000$ and $r$ varying from 2000 to 10000, we test 100 independent instances for each case and record the average CPU time of DPA and NPG. The results are plotted in Fig.1. As can be seen in Fig.1, our DPA is not that sensitive with the sparsity $r$ in terms of the CPU time comparing to NPG. This is reasonable since the computation time for the core projection step tailored for such simplex-constrained SLP in NPG increases with the growth of $r$.

\begin{figure}[ht]
\begin{center}
\includegraphics[width=2.5 in]{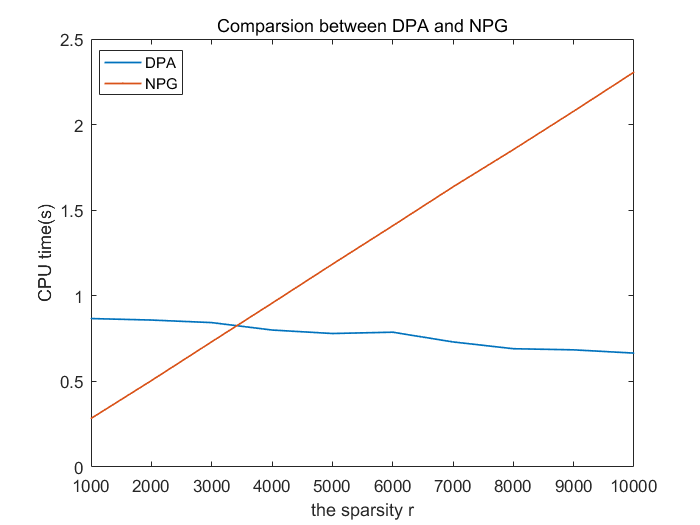}
\end{center}
\caption{The comparison between DPA and NPG with different sparsity constraints}
\end{figure}


%
%
%
%
%
%
%


\section{Conclusions}


In this paper, we have given an explicit formula for the Lagrangian dual of sparse linear programming problems by employing the disjunctive structure of the sparsity constraint involved. The resulting dual problem has been shown to admit the strong duality theorem in the theoretical perspective, and can be efficiently solved by our proposed sPADMM algorithm by fully taking the advantage of the efficient computation for the proximal mapping of the vector Ky-Fan norm function in the computational perspective. Combining with the analysis on optimal index sets of the original problems, we have designed a dual-primal algorithm for pursuing a global solution for a class of sparse linear programming problems. Future research will be focused on those cases in which none optimal index sets can be identified from the dual optimal information. Extensions to sparse nonlinear programming problems deserves further study as well.

\section*{Acknowledgment}
We would like to acknowledge the financial support for the National Natural Science Foundation of China (11771038,11431002,11728101), and the State Key Laboratory of Rail Traffic Control and Safety, Beijing Jiaotong University (RCS2017ZJ001,RCS2016ZZ01).

\bibliographystyle{IEEEtran}
\bibliography{bib_SLP}

\begin{thebibliography}{10}
\providecommand{\url}[1]{#1}
\csname url@samestyle\endcsname
\providecommand{\newblock}{\relax}
\providecommand{\bibinfo}[2]{#2}
\providecommand{\BIBentrySTDinterwordspacing}{\spaceskip=0pt\relax}
\providecommand{\BIBentryALTinterwordstretchfactor}{4}
\providecommand{\BIBentryALTinterwordspacing}{\spaceskip=\fontdimen2\font plus
\BIBentryALTinterwordstretchfactor\fontdimen3\font minus
  \fontdimen4\font\relax}
\providecommand{\BIBforeignlanguage}[2]{{%
\expandafter\ifx\csname l@#1\endcsname\relax
\typeout{** WARNING: IEEEtran.bst: No hyphenation pattern has been}%
\typeout{** loaded for the language `#1'. Using the pattern for}%
\typeout{** the default language instead.}%
\else
\language=\csname l@#1\endcsname
\fi
#2}}
\providecommand{\BIBdecl}{\relax}
\BIBdecl

\bibitem{donoho2006compressed}
D.~L. Donoho, ``Compressed sensing,'' \emph{IEEE Transactions on information
  theory}, vol.~52, no.~4, pp. 1289--1306, 2006.

\bibitem{Liu2009Large}
J.~Liu, J.~Chen, and J.~Ye, ``Large-scale sparse logistic regression,'' in
  \emph{ACM SIGKDD International Conference on Knowledge Discovery and Data
  Mining, Paris, France, June 28 - July}, 2009, pp. 547--556.

\bibitem{Li2006OPTIMAL}
D.~Li, X.~Sun, and J.~Wang, ``Optimal lot solution to cardinality constrained
  mean¨cvariance formulation for portfolio selection,'' \emph{Social Science
  Electronic Publishing}, vol.~16, no.~1, pp. 83--101, 2006.

\bibitem{chen2014sparsity}
A.~I. Chen and S.~C. Graves, ``Sparsity-constrained transportation problem,''
  \emph{arXiv preprint arXiv:1402.2309}, 2014.

\bibitem{babu2010linear}
P.~Babu, K.~Pelckmans, P.~Stoica, and J.~Li, ``Linear systems, sparse
  solutions, and sudoku,'' \emph{IEEE Signal Processing Letters}, vol.~17,
  no.~1, pp. 40--42, 2010.

\bibitem{friedlander2007exact}
M.~P. Friedlander and P.~Tseng, ``Exact regularization of convex programs,''
  \emph{SIAM Journal on Optimization}, vol.~18, no.~4, pp. 1326--1350, 2007.

\bibitem{donoho2005sparse}
D.~L. Donoho and J.~Tanner, ``Sparse nonnegative solution of underdetermined
  linear equations by linear programming,'' \emph{Proceedings of the National
  Academy of Sciences of the United States of America}, vol. 102, no.~27, pp.
  9446--9451, 2005.

\bibitem{shang2014minimal}
M.~Shang, C.~Zhang, and N.~Xiu, ``Minimal zero norm solutions of linear
  complementarity problems,'' \emph{Journal of Optimization Theory and
  Applications}, vol. 163, no.~3, pp. 795--814, 2014.

\bibitem{Chen2016Sparse}
X.~Chen and S.~Xiang, ``Sparse solutions of linear complementarity problems,''
  \emph{Mathematical Programming}, vol. 159, no. 1-2, pp. 539--556, 2016.

\bibitem{d2008first}
A.~d'Aspremont, O.~Banerjee, and L.~El~Ghaoui, ``First-order methods for sparse
  covariance selection,'' \emph{SIAM Journal on Matrix Analysis and
  Applications}, vol.~30, no.~1, pp. 56--66, 2008.

\bibitem{candes2008restricted}
E.~J. Candes, ``The restricted isometry property and its implications for
  compressed sensing,'' \emph{Comptes rendus mathematique}, vol. 346, no. 9-10,
  pp. 589--592, 2008.

\bibitem{Cohen2009COMPRESSED}
A.~Cohen, W.~Dahmen, and R.~Devore, ``Compressed sensing and best k-term
  approximation,'' \emph{Journal of the American Mathematical Society},
  vol.~22, no.~1, pp. 211--231, 2009.

\bibitem{Zhao2013RSP}
Y.~B. Zhao, ``Rsp-based analysis for sparsest and least <formula
  formulatype="inline"> -norm solutions to underdetermined linear systems,''
  \emph{IEEE Transactions on Signal Processing}, vol.~61, no.~22, pp.
  5777--5788, 2013.

\bibitem{beck2015minimization}
A.~Beck and N.~Hallak, ``On the minimization over sparse symmetric sets:
  projections, optimality conditions, and algorithms,'' \emph{Mathematics of
  Operations Research}, vol.~41, no.~1, pp. 196--223, 2015.

\bibitem{xu2016efficient}
F.~Xu, Z.~Lu, and Z.~Xu, ``An efficient optimization approach for a
  cardinality-constrained index tracking problem,'' \emph{Optimization Methods
  and Software}, vol.~31, no.~2, pp. 258--271, 2016.

\bibitem{Panlili2017}
L.~Pan, S.~Zhou, N.~Xiu, and H.~Qi, ``A convergent iterative hard thresholding
  for sparsity and nonnegativity constrained optimization,'' vol.~13, pp.
  325--353, 04 2017.

\bibitem{diewert1974applications}
W.~E. Diewert, ``Applications of duality theory,'' 1974.

\bibitem{Rockafellar1970Convex}
R.~T. Rockafellar, \emph{Convex Analysis}.\hskip 1em plus 0.5em minus
  0.4em\relax Princeton University Press, 1970.

\bibitem{sun2015convergent}
D.~Sun, K.-C. Toh, and L.~Yang, ``A convergent 3-block semiproximal alternating
  direction method of multipliers for conic programming with 4-type
  constraints,'' \emph{SIAM journal on Optimization}, vol.~25, no.~2, pp.
  882--915, 2015.

\bibitem{hardy1952inequalities}
G.~H. Hardy, J.~E. Littlewood, and G.~P{\'o}lya, \emph{Inequalities}.\hskip 1em
  plus 0.5em minus 0.4em\relax Cambridge university press, 1952.

\bibitem{overton1993optimality}
M.~L. Overton and R.~S. Womersley, ``Optimality conditions and duality theory
  for minimizing sums of the largest eigenvalues of symmetric matrices,''
  \emph{Mathematical Programming}, vol.~62, no. 1-3, pp. 321--357, 1993.

\bibitem{watson1993matrix}
G.~Watson, ``On matrix approximation problems with ky fank norms,''
  \emph{Numerical Algorithms}, vol.~5, no.~5, pp. 263--272, 1993.

\bibitem{wu2014moreau}
B.~Wu, C.~Ding, D.~Sun, and K.-C. Toh, ``On the moreau--yosida regularization
  of the vector k-norm related functions,'' \emph{SIAM Journal on
  Optimization}, vol.~24, no.~2, pp. 766--794, 2014.

\bibitem{fazel2013hankel}
M.~Fazel, T.~K. Pong, D.~Sun, and P.~Tseng, ``Hankel matrix rank minimization
  with applications to system identification and realization,'' \emph{SIAM
  Journal on Matrix Analysis and Applications}, vol.~34, no.~3, pp. 946--977,
  2013.

\bibitem{VA}
R.~T. Rockafellar and R.~J.-B. Wets, \emph{Variational Analysis}.\hskip 1em
  plus 0.5em minus 0.4em\relax Springer-Verlag Berlin Heidelherg, 2009.

\end{thebibliography}

\end{document}